\documentclass[10pt]{amsart}
\usepackage{amsmath,amscd}
\usepackage{amsbsy}
\usepackage{amssymb}
\usepackage{amscd,amsthm}
\usepackage[all,cmtip]{xy}
\usepackage[russian,english]{babel}
\usepackage[utf8]{inputenc}

\newtheorem{thm}{Theorem}\numberwithin{thm}{section}

\newtheorem{cor}[thm]{Corollary}

\newtheorem{prob}[thm]{Problem}

\newtheorem*{con2}{Conjecture}

\newtheorem*{exam2}{Example}
\newtheorem*{rema2}{Remark}

\begin{document}
\begin{center}
\huge{A note on some polynomial-factorial diophantine equations}\\[1cm]
\end{center}
\begin{center}

\large{Sa$\mathrm{\check{s}}$a Novakovi$\mathrm{\acute{c}}$}\\[0,5cm]
{\small August 2023}\\[0,5cm]
\end{center}
%\end{center}
\begin{center}
\begin{otherlanguage*}{russian}
	\emph{ЗА APИАНУ.}	
\end{otherlanguage*}
\end{center}
{\small \textbf{Abstract}. 
In 1876 Brocard, and independently in 1913 Ramanujan, asked to find all integer solutions for the equation $n!=x^2-1$. It is conjectured that this equation has only three solutions, but up to now this is an open problem. Overholt observed that a weak form of Szpiro's-conjecture implies that Brocard's equation has finitely many integer solutions. More generally, assuming the ABC-conjecture, Luca showed that equations of the form $n!=P(x)$ where $P(x)\in\mathbb{Z}[x]$ of degree $d\geq 2$ have only finitely many integer solutions with $n>0$. And if $P(x)$ is irreducible, Berend and Harmse proved unconditionally that $P(x)=n!$ has only finitely many integer solutions. In this note we study diophantine equations of the form $g(x_1,...,x_r)=P(x)$ where $P(x)\in\mathbb{Z}[x]$ of degree $d\geq 2$ and $g(x_1,...,x_r)\in \mathbb{Z}[x_1,...,x_r]$ where for the $x_i$ one may also plug in $A^{n}$ or the Bhargava factorial $n!_S$. We want to understand when there are finitely many or infinitely many integer solutions. %Of course, there are plenty of polynomials $g(x_1,...,x_r)$ for which the mentioned diophantine equation has infinitely many integer solutions. Therefore, we restrict ourselves to a certain class of poynomials $g$. We show that for this class of polynomials the ABC-conjecture implies that diophantine equations of the mentioned form have only finitely many integer solutions, generalizing in this way the result of F. Luca, but also a result of M. Ulas. 
Moreover, we study diophantine equations of the form $g(x_1,...,x_r)=f(x,y)$ where $f(x,y)\in\mathbb{Z}[x,y]$ is a homogeneous polynomial of degree $\geq2$. }

%we show that this implies even more, namely that a diophantine equation of the form $bA^n=P(x)$ with $P(x)$ a polynomial in $\mathbb{Z}[x]$ of degree $d\geq 2$ and $b$ and $A$ fixed has only finitely many integer solutions $(n,x)$. %In this way we generalize some previously known result.%We also give examples of w-helices by which we obtain solutions of Markov-type equations. 
\begin{center}
\tableofcontents
\end{center}
\section{Introduction}
Brocard \cite{BR} and independently Ramanujan \cite{RA} asked to find all integer solutions for $n!=x^2-1$. Up to now this is an open problem, known as Brocard's problem. It is believed that the equation has only the three solutions $(x,n)=(5,4), (11,5)$ and $(71,7)$. Overholt \cite{O} observed that a weak form of Szpiro's-conjecture implies that Brocard's equation has finitely many integer solutions. Using the ABC-conjecture Luca \cite{L} proved that diophantine equations of the form $n!=P(x)$ with $P(x)\in \mathbb{Z}[x]$ of degree $d\geq 2$ have only finitely many integer solutions with $n>0$. If $P(x)$ is irreducible, Berend and Harmse \cite{BH1} showed unconditionally that $P(x)=H_n$ has finitely many integer solutions where $H_n$ are highly divisible sequences which also include $n!$. Furthermore, they proved that the same is true for certain reducible polynomials. %In the present work, we consider diophantine equations of the form $bn_1!A^{n_1}\cdots n_r!A_r^{n_r}=P(x)$ with fixed non-zero integer $b$ and fixed positive integers $A_1,...,A_r$ and prove an analogous statement for such equations. %Note that Siegel \cite{S} proved that $bA^n+D=x^2$ has only finitely many integer solutions. 

Without assuming the ABC-conjecture, Berend and Osgood \cite{BO} showed that for arbitrary $P(x)$ of degree $\geq2$ the density of the set of positive integers $n$ for which there exists an integer $x$ such that $P(x)=n!$ is zero. We believe that the arguments presented in \emph{loc.cit}. can also be applied in more general situations, implying, for instance, that the density of the set of positive integers solving $n_1!A^{n_1}\cdots n_r!A_r^{n_r}=P(x)$ is zero. Further progress in this direction was obtained by Bui, Pratt and Zaharescu \cite{BPZ} where the authors give an upper bound on integer solutions $n\leq N$ to $n!=P(x)$. For a detailed overview on results about, for instance, Ramanujan--Nagell type equations $bA^n+D=x^2$ and exponential diophantine equations in general we refer to \cite{BMS}, \cite{SH} and \cite{ST}. Of course, there are several polynomials $P(x)$ for which $P(x)=n!$ is known to have either very few integer solutions or none (see for instance \cite{EO}, \cite{DAB} and \cite{PS}). %But finding all the integer solutions for $x^2-1=n!$ is still an unsolved problem, known as Brocard's problem \cite{BR}. 
Berndt and Galway \cite{BG} showed that the largest value of $n$ in the range $n<10^9$ for which Brocard's equation $x^2-1=n!$ has a positive integer solution is $n=7$. Matson \cite{MA} extended the range to $n<10^{12}$ and Epstein and Glickman \cite{EG} to $n<10^{15}$. %Overholt \cite{O} observed that a weak form of Szpiro's-conjecture actually implies that Brocard's equation has only finitely many integer solutions. 

Starting from Brocard's problem there are also studied variations or generalizations of $x^2-1=n!$ (see for instance \cite{DMU}, \cite{KL}, \cite{MU}, \cite{MUT} \cite{TY}). For instance Ulas \cite{MU} studied, among others, diophantine equations of the form $2^nn!=y^2$ and proved that the Hall conjecture (which is a special case of ABC-conjecture) implies that the equation has only finitely many integer solutions. Note that $2^nn!$ can also be formulated using the notation for the Bhargava factorial $n!_S$. Then $2^nn!=n!_S$, with $S=\{2n+b| n\in \mathbb{Z}\}$. We do not recall the definition of the Bhargava factorial and refer to \cite{BH} or \cite{WT} instead.  

In the present note we are interested in the following problem: let $g(x_1,...,x_r)\in\mathbb{Z}[x_1,...,x_r]$ and $f$ be  polynomials where either $f\in\mathbb{Q}[x]$ or $f\in\mathbb{Q}[x,y]$. Then consider the diophantine equation
\begin{center}
	$g(x_1,...,x_r)=f$.
\end{center}
For any of the $x_i$ in $g$ we may also plug in $A^n$ or $n!_S$. We think that this is a reasonable generalization. Of course, formulated in this generality there are plenty of $g$ and $f$ such that the diophantine equation has infinitely many (positive) integer solutions. Several (exponential) diophantine equations are of the above form. Below we mention only a few examples:
\begin{itemize}
	\item[(i)] \emph{superelliptic equations}.
	\item[(ii)] \emph{Erd\"os-Obl\'ath type equations}: take $g(x_1)=x_1$ and plug in $n!$ and let $g(x,y)=x^p\pm y^p$ or take $g(x_1,x_2)=x_1\pm x_2$ and plug in $n!$ respectively $m!$ and let $f(x)=x^p$.
	\item[(iii)] \emph{Thue-Mahler equation}: take $g(x_1,...,x_r)=x_1\cdot x_2\cdots x_r$ and let $x_i=p_i^{n_i}$.
	\item[(iv)] \emph{Thue-equation}: take $g(x_1,...,x_r)=m$ and let $f(x,y)$ be a homogeneous polynomial of degree $\geq 3$.
	\item[(v)] \emph{Brocard's problem}: take $g(x_1)=x_1$ and plug in $n!$ and let $f(x)=x^2-1$. More generally, let $f(x)$ be any polynomial of degree $\geq 2$ and we get the diophantine equation considered in \cite{L}.
	\item[(vi)] \emph{Ramanujan-Nagell type equation}: take $g(x_1)=bx_1+D$ and plug in $A^n$ for some positive fixed $A$ and let $f(x)=x^2$.
	\item[(vii)] \emph{Fermat equation}: take $g(x_1,x_2)=x_1^n+x_2^n$ and let $f(x)=x^n$.
	\item[(viii)] \emph{generalizations of Brocard's problem}, see \cite{WT}: take $g(x_1)=x_1$ and plug in $n!$ and let $f(x,y)$ be any irreducible binary form of degree $\geq 2$. 
	\item[(ix)] \emph{generalizations of Brocard's problem}, see \cite{DA}: take $g(x_1,x_2)=-x_1^2+x_2$ and plug in $n!$ and let $f(x,y)=x^2+y^2-A$.
\end{itemize}

In the examples (i) to (viii) from above, conditionally or unconditionally there are only finitly many integer solutions. In some situations the exact number of integer solutions is known. So from a structural point of view one can ask for a geometric characterization of the polynomials $g$ and $f$ (or the variety defined by $g-f=0$) such that the diophantine equation $g=f$ has finitely many solutions. % In general, this is a deep problem, tackeld in diophantine geomety. 
For instance, if the hypersurface defined by $g-f$ is of general type (i.e. if the degree is large enough), the Bombieri--Lang conjecture states that the set of rational points is not Zariski dense. So if $g-f$ describes a curve in $\mathbb{P}^2$, Falting's theorem tells us that there are only finitely many rational and therefore integer solutions. Pluging in $n!$ or $n!A^n$ into one of the variables of $g$ still gives only finitely many integer solutions. But often we encounter diophantine equations that are not of general type or that produce varieties of higher dimension where only few structural results concerning the existence and number of rational points are available. At this point, one can also try the so called \emph{modular approach} which was used in Wiles' celebrated proof of Fermat's last theorem. The idea is, roughly speaking, to assign to a diophantine equation a so called \emph{Frey-Hellegouarch curve}. Hopefully, this curve is an elliptic curve with minimal discriminant $\Delta=C\cdot D^p$. Then look at the Galois-representation of the $p$-torsion of the curve and conclude, for instance, by results of Mazur that the representation is irreducible. Now one can apply Ribet's theorem to obtain that the representation gives rise to some newform of level $N$. If one can prove that there are no such newforms of level $N$, one finally finds that the considered diophantine equation has no non-tivial solutions. Besides equations of Fermat-type, this approach was also applied to exponential diophantine equations such as generalized Ramanujan-Nagell-equations. For details we refer the reader to the survey of Siksek \cite{SS} and references therein. 

In the present note we will not use the modular approach. Instead, we  make use of elementary arguments, the ABC or weak form of Szpiro's-conjecture and indirectly of diophantine approximation and methods from algebraic number theory. At some points the proofs can certainly be simplyfied, but I wanted to make the arguments as detailed as necessary. 

We recall the ABC and the weak form of Szpiro's-conjecture which can be found for instance in \cite{LA}. For a non-zero integer $a$, let $N(a)$ be the \emph{algebraic radical}, namely $N(a)=\prod_{p|a}{p}$. Note that 
\begin{eqnarray}
N(a)=\prod_{p|a}{p}\leq \prod_{p\leq a}{p}< 4^a,
\end{eqnarray}
where the last inequality follows from a Chebyshev-type result in elementary prime number theory and is called the Finsler inequality.
\begin{con2}[Weak form of Szpiro's-conjecture]
There exists some constant $s>0$ such that for mutually prime integers $A,B$ and $C$  with $A+B=C$ the inequality 
\begin{eqnarray}
|ABC|<N(ABC)^{s}
\end{eqnarray}
holds.
\end{con2} 
%A generalization of the above conjecture is the following:
\begin{con2}[ABC-conjecture]
For any $\epsilon >0$ there is a constant $K(\epsilon)$ depending only on $\epsilon$ such that whenever $A,B$ and $C$  are three coprime and non-zero integers with $A+B=C$, then 
\begin{eqnarray}
\mathrm{max}\{|A|,|B|,|C|\}<K(\epsilon)N(ABC)^{1+\epsilon}
\end{eqnarray}
holds.
\end{con2} 
The ABC-conjecture from above implies the weak form of Szpiro's-conjecture.\\
%\begin{prop}
%Let $P(x)\in\mathbb{Z}[x]$ a polynomial of degree $d\geq 2$ and $A$ a non-zero integer with $\mathrm{log}(A)>2d+1$. The ABC-conjecture implies that $A^n=P(x)$ has only finitely many integer solutions.
%\end{prop}
%We reprove Siegel's result for $bD^n+E=x^2$ by assuming that Szpiro's-conjecture holds true.
%\begin{prop}
%The weak form of Szpiro's-conjecture implies that $bD^n+E=x^2$ for fixed positive $b,D$ and $E$ has only finitely many integer solutions.
%\end{prop}

In the following let $g(x_1,x_2,...,x_t)=bx_1\cdot x_2\cdots x_t$. In all the variables we plug in $n!$ or $n!A^n$. Formulated in the notation of the Bhargava factorial, we plug in $n!_{\mathbb{Z}}$ or $n!_S$, where $S=\{An+b | n\in\mathbb{Z}\}$ for some fixed positive integer $A$. %The main results of the present note are the following:
Therefore, we want to study diophantine equations of the form
\begin{center}
	$bn_1!\cdots n_r!A_1^{n_1}\cdots A_q^{n_q}=f(x)$
\end{center}
and 
\begin{center}
	$bn_1!\cdots n_r!A_1^{n_1}\cdots A_q^{n_q}=f(x,y)$,
\end{center}
where $q\leq r$, $r>0$ and $f(x)$ and $f(x,y)$ are polynomials with rational coefficients. For a better readibility and to keep it clearer, we prove the results in the present work for $q=r$. However, we would like to note that all results also hold for $q<r$.

\begin{thm}
Fix a non-zero integer $b$ and positive integers $A_1,...,A_r$. If $d>r$, then the equation $bn_1!A_1^{n_1}\cdots n_r!A_r^{n_r}=x^d$ has only finitely many integer solutions. If $d\leq r$, then $bn_1!A_1^{n_1}\cdots n_r!A_r^{n_r}=x^d$ has infinitely many integer solutions, except when $b<0$ and $d$ is even, where there are no solutions.
\end{thm}
\begin{rema2}
	\textnormal{The statement of Theorem 1.1 remains true if $x^d$ is replaced by $ax^d$ with fixed positive rational number $a$.}
\end{rema2}
%Notice that for $s\cdot d=r$ one easily finds infinitely many positive integer solutions for $bn_1!A_1^{n_1}\cdots n_r!A_r^{n_r}=x^d$. If $d<r$ and $\mathrm{gcd}(d,r)=1$, then
\begin{thm}
Let $f(x)\in\mathbb{Q}[x]$ be a polynomial of degree $d\geq 2$ which is not monomial and has at least two distinct roots. Fix a non-zero integer $b$ and positive integers $A_1,...,A_r$. Then the ABC-conjecture implies that $bn_1!A_1^{n_1}\cdots n_r!A_r^{n_r}=f(x)$ has only finitely many integer solutions with $n_i>0$.  
\end{thm}
We want to point out that we can generalize Theorem 1.2 by following its proof. In fact, it is used that $N(bn_1!A_1^{n_1}\cdots n_r!A_r^{n_r})=o(bn_1!A_1^{n_1}\cdots n_r!A_r^{n_r})$ as $(n_1,...n_r)\rightarrow \infty$. We can argue as in the proof of Theorem 1.2 to obtain:

\begin{thm}
	Let $f(x)\in\mathbb{Q}[x]$ be a polynomial of degree $d\geq 2$ which is not monomial and has at least two distinct roots and $F(n_1,...,n_r)$ a function satisfying $N(F(n_1,...n_r))=o(F(n_1,...n_r))$ as $(n_1,...,n_r)\rightarrow \infty$. Then the ABC-conjecture implies that $P(x)=F(n_1,...,n_r)$ has finitely many integer solutions.
\end{thm}
Since $n!|n!_S$, we see that for all primes $p$, the $p$-adic valuation from the definition of $n!_S$ tends to infinity as $n\rightarrow \infty$. Therefore, $N(n!_S)=o(n!_S)$ as $l\rightarrow \infty$. From this, it follows that $N(n_1!_{S_1}\cdots n_r!_{S_r})=o(n_1!_{S_1}\cdots n_r!_{S_r})$ as $(n_1,...,n_r)\rightarrow \infty$. Theorem 1.3 then implies:
\begin{cor}
	Let $f(x)\in\mathbb{Q}[x]$ be a polynomial of degree $d\geq 2$ which is not monomial and has at least two distinct roots. Then the ABC-conjecture implies that  $f(x)=n_1!_{S_1}\cdots n_r!_{S_r}$ has finitely many integer solutions.
\end{cor}
Note that Theorems 1.1, 1.2 and 1.3 generalize the results in \cite{L} and some of the results in \cite{MU} and \cite{WT}. 
Dabrowski \cite{DA} asked whether equations of the form $n!+A=x^2+y^2$ have finitely many positive integer solutions and Ulas \cite{MU} studied equations of the form $x^2-A=n!!$ %Furthermore, Berend and Hermse \cite{BH} studied equations of the form $x^r(x+1)=n!$ and proved that for $x\geq 4$ there are only finitely many postivie integer solutions. 
%In the present note we are also interested in equations of the form $n!!=x^2y\pm y^2x$ (see Example on page 10), 
where $n!!$ denotes the double factorial. Other equations involving the double factorial have also been studied in \cite{MU}. As mentioned above, Luca \cite{L} proved uncoditionally that $x^d=n!$ has only finitely many solutions for $d\geq 2$. Takeda \cite{WT} considered, more generally, equations of the form $f(x,y)=n!$ where $f(x,y)$ is an arbitrary  binary form of degree $\geq2$. He proved, among others, that if $f(x,y)$ is irreducible of degree $\geq 2$, then there are only finitely many $n$ such that $n!$ is represented by $f(x,y)$. With Thue's theorem one has finitely many solutions if the degree of $f(x,y)$ is at least three. In \emph{loc.cit.} it is furthermore observed that the same holds for certain reducible polynomials $f(x,y)$. So in view of this fact and, for instance, the work of Erd\"os and Obl\'ath \cite{EO}, it is reasonable to ask, for example, whether $f(x,y)=n!_S$ or $f(x,y)=n!!$ have only finitely many integer solutions when $f(x,y)$ is a (homogeneous) polynomial of degree at least two. 

Theorems 1.5, 1.7, 1.8 and 1.12 make a small step in this direction.
\begin{thm}
	Let $f(x,y)=a_dx^d+a_{d-1}x^{d-1}y+\cdots + a_{1}xy^{d-1}+a_0y^d\in \mathbb{Z}[x,y]$ be an irreducible homogeneous polynomial. Fix some non-zero integer $b$ and positive integers $A_1,...,A_r$. If $d>r$, then there are finitely many $(n_1,...,n_r)$ such that $	bn_1!A_1^{n_1}\cdots n_r!A_r^{n_r}$ is represented by $f(x,y)$. If $r\geq 2$, then the diophantine equation
\begin{eqnarray*}
	bn_1!A_1^{n_1}\cdots n_r!A_r^{n_r}=f(x,y)
\end{eqnarray*}
has finitely many integer solutions with $n_i>0$. 	
\end{thm}

\begin{cor}
	For any irreducible $f(x)\in\mathbb{Z}[x]$ with degree $d>r$, the equation $f(x)=bn_1!A_1^{n_1}\cdots n_r!A_r^{n_r}$ has only finitely many integer solutions.
	\end{cor}
Theorems 1.7, 1.8 and 1.12 deal with the case of some reducible polynomials. 
\begin{thm}
	Let $f(x,y)\in\mathbb{Z}[x,y]$ be a polynomial with factorization 
	\begin{center}
		$f(x,y)=f_1(x,y)^{e_1}\cdots f_u(x,y)^{e_u}$,
	\end{center}
	where $f_i(x,y)$ are irreducible homogeneous polynomials of degree $d_i$. Fix some non-zero integer $b$ and positive integers $A_1,...,A_r$. Now if $d_i\geq2$ and $d_1e_1+\cdots +d_ue_u>r$ or if $\mathrm{min}\{d_1e_1,...,d_ue_u\}>r$, then there are only finitely many $(n_1,...,n_r)$ such that $=bn_1!A_1^{n_1}\cdots n_r!A_r^{n_r}$ is represented by $f(x,y)$. If $r\geq 2$, then the equation 
	\begin{eqnarray*}
		bn_1!A_1^{n_1}\cdots n_r!A_r^{n_r}=f(x,y)
	\end{eqnarray*}
	has only finitely many integer solutions.
\end{thm}
Theorem 1.7 treats equations such as $f(x,y)=(xy)^4(x-y)^3=n!m!$ but can not be applied to $xy(x\pm y)=3^nn!$. Theorems 1.8 and 1.12 below  try to include cases where Theorem 1.7 can not be applied (see the examples on page 13 and 14).
\begin{thm}
Let $P(x), Q(y)$ be polynomials with rational coefficients such that at least one of the two has at least two distinct roots. % and such that for any pair $(\alpha,\beta)$ with $\mathrm{gcd}(\alpha,\beta)=1$ one has $\mathrm{gcd}(P(\alpha),Q(\beta)))=1$. 
Set $f(x,y):=P(x)^2Q(y)+Q(y)^2P(x)$ and fix some non-zero integer $b$ and positive integers $A_1,...,A_r$. Then the ABC-conjecture implies that 
\begin{eqnarray*}
bn_1!A_1^{n_1}\cdots n_r!A_r^{n_r}=f(x,y)
\end{eqnarray*}
has only finitely many positive integer solutions $(n_1,...,n_r,x,y)$ with $n_i>0$.  %$\mathrm{gcd}(x,y)=1$ ?. 
\end{thm}
We want to stress that the assumption of %coprime $x$ and $y$ and 
the existence of at least two distinct roots is necessary to conclude that there are only finitely many positive integer solutions. Consider for instance $P(x)=x$ and $Q(y)=y$ and the equation
\begin{center}
$P(x)^2Q(y)+Q(y)^2P(x)=2\cdot n!\cdot m!\cdot l!$	
\end{center}
Choosing $x=y$, we see that 
\begin{center}
	$2x^3=2\cdot n!\cdot m!\cdot l!$
\end{center}
has infinitely many solutions according to Theorem 1.1. In this context we make the observation that Theorem 1.1 has the following consequence.
\begin{cor}
	Let $f(x,y)=a_dx^d+a_{d-1}x^{d-1}y+\cdots + a_{1}xy^{d-1}+a_0y^d$ be a homogeneous polynomial of degree $d$. Assume $a_0+\cdots +a_d>0$ or $a_0>0$ or $a_d>0$. Furthermore, let $b>0$ and $d\leq r$. Then the diophantine equation 
\begin{eqnarray*}
	bn_1!A_1^{n_1}\cdots n_r!A_r^{n_r}=f(x,y).
\end{eqnarray*}	
has infinitely many integer solutions with $n_i>0$. 
\end{cor}
There are certainly more cases where the equation of Corollary 1.8 has infinitely many solutions. Our aim was just to point out that in case $d\leq r$ there are usually infinitely many solutions. Notice that Theorem 1.1 also implies that the equation in Corollary 1.8 has only finitely many integer solutions with $x=y$ if $d>r$. We do not know whether there are only finitely many solutions with $x\neq y$ and $\mathrm{gcd}(x,y)\neq 1$ in case $d>r$. 
\begin{prob}
	\textnormal{Let $f(x,y)$ be a homogeneous polynomial of degree $d$ with rational coefficients. Fix a non-zero integer $b$ and some positive integers $A_1,...,A_q$ and assume $d>r$. Furthermore, let $S_1,...,S_r$ be infinite subsets of $\mathbb{Z}$. Does the equation $bn_1!_{S_1}\cdots n_r!_{S_r}\cdot A_1^{m_1}\cdots A_q^{m_q}=f(x,y)$ have only finitely many integers solutions $(n_1,...,n_r,m_1,...,m_q,x,y)$? Are there infinitely many solutions if $d\leq r$?}
\end{prob}
A related problem is the following:
\begin{prob}
	\textnormal{Let $f(x,y)$ be a homogeneous polynomial of degree $d$ with rational coefficients. Let $S_1,...,S_r$ be infinite subsets of $\mathbb{Z}$ and fix some non-zero integer $b$. Give some characterization for when $f(x,y)=bl_1!_{S_1}\cdots l_r!_{S_r}$ has finitely many or infinitely many integer solutions? }
\end{prob}
%\begin{exam2}
%Let $x^3d\pm y^3d=n!!$. Then for any fixed positive integer $d$ the ABC-conjecture implies that there are only finitely many positive integer solutions.	
%\end{exam2}

Another result is this direction is obtained if $A_1,...,A_q$ in Problem 1.9 is a set of distinct primes. Using the Thue-Mahler Theorem, we can show the following:
\begin{thm}
Let $f(x,y)=x^{2s}y^s\pm y^{2s}x^s$ where $s>0$ is an arbitrary but fixed positive integer %such that for any pair $(\alpha,\beta)$ with $\mathrm{gcd}(\alpha,\beta)=1$ one has $\mathrm{gcd}(P(\alpha),Q(\beta)))=1$ 
and $p_1,...,p_q$ a set of primes with $p_1<...<p_q$.  Fix some non-zero integer $b$ and assume $f(x,y)$ has at least three pairwise non-proportional linear factors over $\mathbb{C}$. Then the weak form of Szpiro's-conjecture implies that 
\begin{eqnarray*}
bp_1^{z_1}\cdots p_q^{z_q}\cdot n_1!A_1^{n_1}\cdots n_l!A_l^{n_l}=f(x,y)
\end{eqnarray*}
has only finitely many integer solutions $(z_1,...,z_q,n_1,...,n_l,x,y)$ with $z_i\geq 0$ and $n_j>0$ and $\mathrm{gcd}(x,y)=1$.
\end{thm}
%Finnaly we summarize the results in the following table:\\

%\begin{tabular}{ r|c|c| }
%	\multicolumn{1}{r}{}
%	&  \multicolumn{1}{c}{$f(x,y)$}
%	& \multicolumn{1}{c}{$f(x)$} \\
%	\cline{2-3}
%	 & homogeneous and irreducible with degree $d>r$ & irreducible of degree $d>r$ or of the form $\frac{a}{b}x^d$\\
%	 &finitely many solutions & finitely many solutions \\
%	\cline{2-3}
%	 & Book & Tutor \\
%	\cline{2-3}
%	&a &b\\
%	\cline{2-3}
%\end{tabular}\\

\noindent
{\small \textbf{Acknowledgement}. I am grateful to Florian Luca for very helpful explanations and comments.}

\section{Proof of Theorem 1.1}
We give the proof only for the cases $r=1,2$, since the arguments for arbitrary $r>2$ are analogous. So let $r=1$. We consider the equation
\begin{center}
$bn!A^n=x^d$.
\end{center}
If $n>2\mathrm{max}\{A,|b|\}$, then there is a prime number $p$ in the interval $(n/2,n)$ wich is larger than $\mathrm{max}\{A,|b|\}$. The prime $p$ will appear with exponent one in $bn!A^n$. Since $d>1$, the number $bn!A^n$ cannot be a perfect power. 

Now let $r=2$. We consider
\begin{eqnarray}
b n! A^n m! B^m=x^d.
\end{eqnarray} 
If $n>2\mathrm{max}\{A,B,|b|\}$, then there is a prime number $p$ in the interval $(n/2,n)$ which is larger than $\mathrm{max}\{A,B,|b|\}$. There are three cases to consider. 
\begin{itemize}
\item[1)] $n>m$. In this case we see that the prime $p$ will appear with exponent at most two in the product $b n! A^n m! B^m$. Since $d\neq 2$, the product $b n! A^n m! B^m$ cannot be of the form $x^d$. Therefore, there are no integer solutions if $n>2\mathrm{max}\{A,B,|b|\}$ and $n>m$.
\item[2)] $n<m$. In this case $m>2\mathrm{max}\{A,B,|b|\}$. Then there is a prime number $p$ in the interval $(m/2,m)$ which is larger than $\mathrm{max}\{A,B,|b|\}$. Again, this prime $p$ will appear with exponent at most two in the product $b n! A^n m! B^m$. Since $d\neq 2$, the product $b n! A^n m! B^m$ cannot be of the form $x^d$. Hence, there are no integer solutions if $n>2\mathrm{max}\{A,B,|b|\}$ and $n<m$.
\item[3)] $n=m$. In this case equation (4) becomes
\begin{eqnarray}
b(n!)^2(AB)^n=x^d.
\end{eqnarray}
As $n>2\mathrm{max}\{A,B,|b|\}$, the prime number $p$ in the interval $(n/2,n)$ will appear with exponent two in the product $b(n!)^2(AB)^n$. Since $d\neq 2$, the product $b(n!)^2(AB)^n$ cannot be of the form $x^d$. Therefore, there are no integer solutions if $n>2\mathrm{max}\{A,B,|b|\}$ and $n=m$.
\end{itemize}
So we are left with $n\leq 2\mathrm{max}\{A,B,|b|\}$. Again, we consider the following cases.
\begin{itemize}
\item[1)] $m<n$. In this case there can be only finitely many integer solutions $(n,m,z)$ satisfying (4).
\item[2)] $n=m$. In this case there can be only finitely many integer solutions $(n,n,z)$ satisfying (4).
\item[3)] $n<m$. In this case, either we must have $n<m\leq2\mathrm{max}\{A,B,b\}$ or $n\leq 2\mathrm{max}\{A,B,b\}<m$. Clearly, if $n<m\leq2\mathrm{max}\{A,B,|b|\}$, there can be only finitely many integer solutions $(n,m,z)$ satisfying (4). Now if $n\leq 2\mathrm{max}\{A,B,|b|\}<m$, we conclude from 2) from above that there is a prime number $p$ in the interval $(m/2,m)$ which is larger than $\mathrm{max}\{A,B,|b|\}$. This prime $p$ will appear with exponent at most two in the product $b n! A^n m! B^m$. Since $d\neq 2$, the product $b n! A^n m! B^m$ cannot be of the form $x^d$. 
\end{itemize}
Summarizing, we see that the equation (4) can have only finitely many integer solutions. This completes the first part of the proof.

Now let us consider the case $d\leq r$. Obviously, there are infinitely many integer solutions for $x=bn_1!A_1^{n_1}n_2!A_2^{n_2}\cdots n_r!A_r^{n_r}$. Therefore, we assume $d\geq2$. Notice that in case $b<0$ and $d$ is even the equation has no solution. So we consider the equation
\begin{center}
	$x^d=bn_1!A_1^{n_1}n_2!A_2^{n_2}\cdots n_r!A_r^{n_r}$,
\end{center}
where $b>0$ and $d$ arbitrary or $b<0$ and $d$ odd. %We only treat the case where $\mathrm{gcd}(r,d)=1$ since the other cases follow easily from that. 
Since $d\leq r$, we can rewrite the equation as
\begin{center}
		$x^d=bn_1!A_1^{n_1}n_2!A_2^{n_2}\cdots n_d!A_d^{n_d}\cdot(n_{d+1}!A_{d+1}^{n_{d+1}}\cdots n_r!A_r^{n_r})$.
\end{center}
Now we set $n_{d+1}=n_{d+2}=\cdots =n_r=1$ and $n_1=n_2=\cdots n_{d-1}=m$ and $n_d=m+1$. Then the equation becomes
\begin{center}
	$x^d=(A_1\cdots A_d)^m\cdot (m!)^d\cdot bA_d\cdot A_{d+1}\cdots A_r\cdot(m+1)$.
\end{center}
We rewrite again:
\begin{center}
		$x^d=(A_1\cdots A_d)^{(m-(d-1))}\cdot (m!)^d\cdot b\cdot(A_1\cdots A_d)^{(d-1)}\cdot A_d\cdot A_{d+1}\cdots A_r\cdot(m+1)$.
\end{center}
Now we want to choose $m$ such that $m-(d-1)=m-d+1=ds$. This is equivalent to $m+1=d(s+1)$. If $b>0$, we set
\begin{center}
	$R:=b\cdot(A_1\cdots A_d)^{(d-1)}\cdot A_d\cdot A_{d+1}\cdots A_r$
\end{center}
Then the above equation becomes 
\begin{center}
		$x^d=(A_1\cdots A_d)^{(m-(d-1))}\cdot (m!)^d\cdot Rd\cdot(s+1)$.
\end{center}
Now we can set $s=(Rd)^{td-1}-1$ where $t>0$ is any positive integer and see that $m-(d-1)=d\cdot((Rd)^{td-1}-1)$ is a multiple of $d$. Notice that $d\geq2$ by assumption and hence $s\geq 1$. Our diophantine equation becomes
\begin{center}
	$x^d=((A_1\cdots A_d)^{((Rd)^{td-1}-1)})^d\cdot ((d(Rd)^{td-1}-1)!)^d\cdot ((Rd)^t))^d$.
\end{center}
This shows that we can find infinitely many interger solutions $(x,n_1,...,n_r)$ with $n_i>0$. If $b<0$ and $d$ odd, we set 
\begin{center}
	$R':=|b|\cdot(A_1\cdots A_d)^{(d-1)}\cdot A_d\cdot A_{d+1}\cdots A_r$. 
\end{center}
Then $R=(-1)R'$ and the diophantine equation becomes
\begin{center}
	$x^d=((A_1\cdots A_d)^{((Rd)^{td-1}-1)})^d\cdot ((d(Rd)^{td-1}-1)!)^d\cdot ((Rd)^t))^d\cdot (-1)^{td}$.
\end{center}
Choosing $t$ odd, we conclude that $td$ is odd. Thus $td-1$ is even. Since $d\geq2$ it follows that $(Rd)^{td-1}-1>0$. This shows that there are also infinitely many integer solutions in this case. Note that the solutions are constructed only using the fixed integers $b,A_1,...,A_r$ and the given degree $d$. This completes the proof.

\section{Proof of Theorem 1.2} 
Multiplying the equation $bn_1!A_1^{n_1}\cdots n_r!A_r^{n_r}=f(x)$ by a certain integer, we may assume that $f(x)$ is a polynomial with integer coefficients. So without loss of generality, we assume 
\begin{eqnarray*}
f(x)=a_0x^d+a_1x^{d-1}+...+a_d
\end{eqnarray*}
with $a_i\in \mathbb{Z}$. Now multiply the equation $bn_1!A_1^{n_1}\cdots n_r!A_r^{n_r}=f(x)$ by $d^da_0^{d-1}$. We obtain
\begin{eqnarray*}
y^d+b_1y^{d-1}+...+b_d=c(n_1!A_1^{n_1}\cdots n_r!A_r^{n_r})
\end{eqnarray*}
for a constant $c$, where $c=bd^da_0^{d-1}$ and $y:=a_0dx$. Notice that $b_i=d^ia_ia_0^{i-1}$ so that we can make the change of variable $z:=y+\frac{b_1}{d}$. Since we are assuming that $f(x)$ has at least two distinct roots, the change of variable produces a polynomial that does not have a monomial of degree $d-1$. Therefore we get the following equation
\begin{eqnarray}
z^d+c_2d^{d-2}+...+c_d=c(n_1!A_1^{n_1}\cdots n_r!A_r^{n_r}).
\end{eqnarray}
Notice that $c_i$ are integer coefficients wich can be computed in terms of $a_i$ and $d$. Now let $Q(X)=X^d+c_2X^{d-2}+...+c_d$ and notice that when $|z|$ is large one has
\begin{eqnarray}
\frac{|z|^d}{2}<|Q(z)|<2|z|^d.
\end{eqnarray}
For the rest of the proof we denote by $C_1,C_2,...$ computable positive constants depending on the coefficients $a_i$ and eventually on some small $\epsilon >0$ which comes into play later by applying the ABC-conjecture.
 
Whenever $(n_1,...,n_r,z)$ is a solution to $n_1!A_1^{n_1}\cdots n_r!A_r^{n_r}=f(x)$ we conclude from (6) and (7) that there exist constants $C_1$ and $C_2$ such that
\begin{eqnarray}
|d\cdot\mathrm{log}|z|-\mathrm{log}(n_1!A_1^{n_1}\cdots n_r!A_r^{n_r})|<C_1,
\end{eqnarray}
for $|z|>C_2$ (see \cite{L} equation (10)). %For technical reasons, we assume that $C_2$ is large enough with respect to $C_1$. 
Now let $R(X)\in \mathbb{Z}[X]$ be such that $Q(X)=X^d+R(X)$. Since $f(x)$ is not monomial and has at least two distinct roots, $R(X)$ can be assumed to be non-zero, let $j\leq d$ be the largest integer with $c_j\neq 0$. We rewrite (6) as
\begin{eqnarray*}
z^j+c_2z^{j-2}+...+c_j=\frac{c(n_1!A_1^{n_1}\cdots n_r!A_r^{n_r})}{z^{d-j}}.
\end{eqnarray*}
Let $R_1(X)$ be the polynomial 
\begin{eqnarray*}
R_1(X):= \frac{R(X)}{X^{d-j}}=c_2X^{j-2}+\cdots +c_j.
\end{eqnarray*}
It is shown in \cite{L} there are constants $C_3$ and $C_4\geq C_2$ such that
\begin{eqnarray*}
0<|R_1(z)|< C_3|z|^{j-2},
\end{eqnarray*}
for $|z|> C_4$. So we have 
\begin{center}
$z^j+R_1(z)=\frac{c(n_1!A_1^{n_1}\cdots n_r!A_r^{n_r})}{z^{d-j}}$. 
\end{center}
For $D=\mathrm{gcd}(z^j, R_1(z))$ we have
\begin{eqnarray*}
\frac{z^j}{D}+\frac{R_1(z)}{D}=\frac{c(n_1!A_1^{n_1}\cdots n_r!A_r^{n_r})}{z^{d-j}D}
\end{eqnarray*}
Applying the ABC-conjecture to $A=\frac{z^j}{D}$, $B=\frac{R_1(z)}{D}$ and $C=\frac{c(n_1!A_1^{n_1}\cdots n_r!A_r^{n_r})}{z^{d-j}D}$, we find
\begin{eqnarray}
\frac{|z|^j}{D}< C_5N(\frac{z^jR_1(z)c(n_1!A_1^{n_1}\cdots n_r!A_r^{n_r})}{D^3})^{1+\epsilon},
\end{eqnarray}
where $C_5$ depends only on $\epsilon$. It is shown in \cite{L}, p.272 that 
\begin{eqnarray}
N(\frac{|z|^j}{D})\leq |z|,\\
N(\frac{R_1(z)}{D})<\frac{C_3|z|^{j-2}}{D}.
\end{eqnarray}
Moreover, we have
\begin{eqnarray*}
N(\frac{c(n_1!A_1^{n_1}\cdots n_r!A_r^{n_r})}{z^{d-j}D})\leq N(c)N(n_1!A_1^{n_1}\cdots n_r!A_r^{n_r})\leq N(c)N(n_1!A_1^{n_1})\cdots N(n_r!A_r^{n_r}).
\end{eqnarray*}
%Now there is a $M\in \mathbb{N}$ such that $A_i^{n_i}\leq n_i!$ for all $n_i>M$. 
This gives 
\begin{eqnarray*}
N(\frac{c(n_1!A_1^{n_1}\cdots n_r!A_r^{n_r})}{z^{d-j}D})\leq N(c)N(n_1!A_1^{n_1})\cdots N(n_r!A_r^{n_r})\leq C_6N(n_1!)\cdots N(n_r!),
\end{eqnarray*} 
where $C_6=N(c)N(A_1)\cdots N(A_r)$. 
From (1) it follows
\begin{eqnarray}
N(\frac{c(n_1!A_1^{n_1}\cdots n_r!A_r^{n_r})}{z^{d-j}D})<C_64^{n_1}\cdots 4^{n_r}=C_64^{(n_1+\cdots +n_r)}
\end{eqnarray}
and from (10), (11) and (12) we get
\begin{eqnarray}
N(\frac{|z|^j}{D})N(\frac{R_1(z)}{D})N(\frac{c(n_1!A_1^{n_1}\cdots n_r!A_r^{n_r})}{z^{d-j}D})<\frac{C_3C_6|z|^{j-1}4^{(n_1+\cdots +n_r)}}{D}
\end{eqnarray}
From inequalities (9) and (13), we obtain
\begin{eqnarray}
\frac{|z|^j}{D}<C_7\bigl(\frac{|z|^{j-1}4^{(n_1+\cdots +n_r)}}{D}\bigr)^{(1+\epsilon)}
\end{eqnarray}

If we choose $\epsilon =\frac{1}{2d}\leq \frac{1}{2j}$, inequality (14) implies that
\begin{eqnarray*}
|z|^{1/2}<|z|^{1+\epsilon -\epsilon j}< C_84^{(n_1+...+n_r)(1+\epsilon)},
\end{eqnarray*}
or simply 
\begin{eqnarray*}
\mathrm{log}|z|<C_9n_1+...+C_9n_r+C_{10}.
\end{eqnarray*}
Thus
\begin{eqnarray*}
d\cdot\mathrm{log}|z|<C_{11}n_1+...+C_{11}n_r+C_{12}.
\end{eqnarray*}
This gives
\begin{eqnarray}
\mathrm{log}(n_1!A_1^{n_1}\cdots n_r!A_r^{n_r})<C_1+d\cdot \mathrm{log}|z|< C_{11}n_1+...+C_{11}n_r+C_{13}.
\end{eqnarray}
We can simplify (15) and finally obtain
\begin{center}
$\mathrm{log}(n_1!)+\mathrm{log}(n_2!)+\cdots +\mathrm{log}(n_r!)<C_{14}n_1+C_{14}n_2+\cdots C_{14}n_r+C_{13}$.	
	\end{center}
Now we can conclude that only finitely many $(n_1,...,n_r)$ satisfy (15). We give the argument for $r=2$. So lets consider an inequality of the form
\begin{center}
	$\mathrm{log}(n!)+\mathrm{log}(m!)<A'n+B'm+C'$
\end{center}
where $A',B'$ and $C'$ are positive constant intergers. Assume there are infinitely many pairs $(n,m)$ of natural numbers satisfying the inequality. There are three cases: 
\begin{itemize}
	\item[1)] infinitely many $n$ and finitely many $m$: let $s$ denote the maximum of these $m$ and $t$ the minimum. Then we have
	\begin{center}
		$\mathrm{log}(n!)+\mathrm{log}(t!)<A'n+B's+C'$
	\end{center}
and therefore
\begin{center}
	$\mathrm{log}(n!)<A'n+E'$
\end{center}
where $E'=\mathrm{log}(t!)+B's+C'$ is a constant. With Stirling's formula for approximating $n!$, we find that there are only finitely many $n$ satisfying 
\begin{center}
	$\mathrm{log}(n!)<A'n+E'$.
\end{center}
This gives a contradition. 
	\item[2)] infinitely many $m$ and finitely many $n$: reverse the role of $n$ and $m$. 
	\item[3)] infinitely many $n$ and infinitely many $m$: the argument is similar. At some point $\mathrm{log}(n!)$ exceeds $A'n$ and $\mathrm{log}(m!)$ exceeds $B'm+C'$. Therefore infinitely many $n$ and infinitely many $m$ is impossible. 
\end{itemize}

Since there are only finitely many $(n_1,...,n_r)$ satisfying  (15), we finally conclude from (8) that $|z|< C_{15}$. This completes the proof.
\section{Proof of Theorems 1.5 and 1.7}
To prove the statement of Theorem 1.5 we have to imitate the proof of Theorem 4.1 in \cite{WT}. We follow the notation of \cite{WT}. Let $f(x,y)=a_dx^d+a_{d-1}x^{d-1}y+\cdots + a_{d-1}xy^{1}+a_0y^d$ be an irreducible polynomial and let $K_F$ be the splitting field of $f(x,1)$. Denote by $C_F$ the set of conjugacy classes of the Galois group $G_F=\mathrm{Gal}(K_F/\mathbb{Q})$ whose cycle type $[h_1,...,h_s]$ satisfies $h_i\geq 2$. For a cycle $\sigma$, the cycle type is defined as the ascending ordered list $[h_1,...,h_s]$ of the sizes of cycles in the cycle decomposition of $\sigma$. For further details we refer to Chapters 2, 3 and 4 in \cite{WT}. Of particular interest are the proofs of Lemma 2.1, Lemma 3.1, Theorem 3.6 and Theorem 4.1. We proceed with our proof. Since $d>r\geq 1$, we conclude from \cite{WT}, Lemma 2.1 that $C_F\neq \emptyset$. Now let $C\in C_F$ be a fixed conjugacy class of $G_F$. We say that a prime $p$ corresponds to $C$ if the Frobenius map $(p,K_F/\mathbb{Q})$ belongs to $C$ (see \cite{WT}, chapter 2 for details). Let $g=\mathrm{gcd}(a_d,...,a_0)$ and $N=gp_1\cdots p_uq_1^{l_1}\cdots q_v^{l_v}$ where $q_i$ are primes corresponding to a conjugacy class in $C_F$ satisfying $\mathrm{gcd}(q,a\Delta_{mod})=1$ where $a\in\{a_d,a_0\}$ and $p_j$ are the other primes (see \cite{WT}, Lemma 3.1 for details).  Here $\Delta_{mod}$ denotes a certain modified discriminant of $f(x,1)$ and is defined by 
\begin{center}
	$\Delta_{mod}=\frac{\Delta_{f(x,1)}}{\mathrm{gcd}(a_n,...,a_0)^{2n-2}}$.
\end{center}
The assumption $d\geq 2$ and \cite{WT}, Lemma 3.1 then imply that if $N$ is represented by $f(x,y)$ and $q|N$ for a prime $q$ corresponding to $C$ satisfying $\mathrm{gcd}(q,a\Delta_{mod})=1$, then $N$ is divisible by $q^d$ at least. Without loss of generality, let $n_r\geq n_{r-1}\geq \cdots \geq n_1$ and assume $n_r$ is big enough, say $n_r>2\mathrm{max}\{|b|,A_1,...A_r\}$. Since the second smallest positive integer divisible by $q$ is $2q$, there is no solution for $bn_1!A_1^{n_1}\cdots n_r!A_r^{n_r}=f(x,y)$ if $q<n_r<2q$. Indeed, since $q<n_r<2q$ we see that $2q<2n_r$ and hence $q<n_r<2q<2n_r$. This implies $\frac{n_r}{2}<q<n_r$. Since  $n_r>2\mathrm{max}\{|b|,A_1,...A_r\}$ and since $d>r$, we see that $q^d$ does not divide $bn_1!A_1^{n_1}\cdots n_r!A_r^{n_r}$. Now apply \cite{WT}, Theorem 3.6 (as in the proof of Theorem 4.1 in \emph{loc.cit.}) to conclude that there exists a prime $q'$ corresponding to $C$ with $q'\in (q,2q)$ and satisfying $\mathrm{gcd}(q',a\Delta_{mod})=1$. Therefore, by the same argument as before, if $q'<n_r<2q'$ there are no integer solutions for  $bn_1!A_1^{n_1}\cdots n_r!A_r^{n_r}=f(x,y)$. By an induction argument we conclude that whenever $n_1>q$ there are no integer solutions. This shows that there are only finitely many $(n_1,...,n_r)$ such that $bn_1!A_1^{n_1}\cdots n_r!A_r^{n_r}$ is represented by $f(x,y)$. Now if $r\geq 2$, we have $d\geq 3$ and we can use Thue's theorem to conclude that there are indeed only finitely many integer solutions.\\
The proof of Theorem 1.7 follows exactly the lines of the proofs of Theorems 4.3 and 4.7 in \cite{WT}. 

If all $d_i\geq 2$, then proceed as follows: let $K_{F_j}$ be the splitting field of $f_j(x,1)$ and denote by $C_{F_j}$ be the set of conjugacy classes whose cycle type has sizes $\geq 2$. Now let $q$ be a prime corresponding to a conjugacy class $C\in \cap_{j=1}^uC_{F_j}$ satisfying $\mathrm{gcd}(q,a\Delta_{mod})=1$.Assume $N$ is represented by $f(x,y)$. Therefore, there are $x_0,y_0$ such that $q|f(x_0,y_0)$. Then there exists a polynomial $f_j(x,y)=a_{j,d_j}x^{d_j}+\cdots + a_{j,0}y^{d_j}$ such that $q|f_j(x_0,y_0)$. Now $\mathrm{gcd}(q,a\Delta_{mod})=1$ with  $a\in\{a_d,a_0\}$ implies $\mathrm{gcd}(q,a(j)\Delta_{j, mod})=1$ where $a(j)\in\{a_{j,d_j},a_{j,0}\}$. Here $\Delta_{j, mod}$ denotes the modified discriminant of $f_j(x,1)$. It follows from \cite{WT}, Lemma 3.1 that $q^d|f(x,y)$. As in the proof of Theorem 1.5 above, we assume without loss of generality that $n_r\geq \cdots \geq n_1$ and choose $n_r$ big enough. Again, one concludes the statement from an inductions argument. The rest of the proof is the same as in the proof of Theorem 1.5. Now let $\mathrm{min}\{d_1e_1,...,d_ue_u\}=d_{i_0}e_{i_0}$ and assume $d_{i_0}=1$ and $e_{i_0}>1$. Then one can argue as in \cite{WT}, Theorem 4.7 to conclude that, let's say, $q$ divides $f(x,y)$ at least $e_{i_0}>1$ times. Again, the rest of the proof is as in the proof of Theorem 1.5 above.

\section{Proof of Theorem 1.8}
For a solution $(n_1,...,n_r,\alpha,\beta)$ we consider $\frac{P(\alpha)}{D}$ and $\frac{Q(\beta)}{D}$ where $D:=D(\alpha,\beta)=\mathrm{gcd}(P(\alpha),Q(\beta))$. Applying the weak form of Szpiro's conjecture to $A=\frac{P(\alpha)}{D}$, $B=\frac{Q(\beta)}{D}$ and $C=\frac{P(\alpha)+Q(\beta)}{D}$, we find
%we have $\mathrm{gcd}(P(x),Q(y))=1$ by assumption. This implies that $P(x), Q(y)$ and $P(x)+Q(y)$ are mutually coprime. Applying the weak form of Szpiro's conjecture to $A=P(x), B=Q(y)$ and $C=P(x)+Q(y)$ we obtain
\begin{eqnarray*}
\Bigl|\frac{bn_1!A_1^{n_1}\cdots n_r!A_r^{n_r}}{D^3}\Bigr|=\Bigl|\frac{f(\alpha, \beta)}{D^3})\Bigr|<H\cdot N(n_1!...n_r!)^s< H\cdot (4^s)^{n_1}\cdots (4^s)^{n_r}.
\end{eqnarray*}
Here $H=N(A_1\cdots A_r)$ is a constant. Now Stirling's formula (see \cite{O}) yields
\begin{eqnarray*}
4^r(\frac{n_1}{e})^{n_1}\cdots (\frac{n_r}{e})^{n_r}\leq n_1!\cdots n_r!
\end{eqnarray*}
and since
\begin{eqnarray*}
	n_1!\cdots n_r!\leq \Bigl|bn_1!A_1^{n_1}\cdots n_r!A_r^{n_r}\Bigr|
\end{eqnarray*}
we obtain
\begin{eqnarray*}
(\frac{n_1}{e})^{n_1}\cdots (\frac{n_r}{e})^{n_r}<D^3\cdot H\cdot (4^{s\cdot n_1-1})\cdots (4^{s\cdot n_r-1}).
\end{eqnarray*}
Arguments similar to the arguments in the proof of Theorem 1.2 yield that this inequality is satisfied only for finitely many $(n_1,...,n_r)$. Now the ABC-conjecture applied to $A=\frac{P(\alpha)}{D}$, $B=\frac{Q(\beta)}{D}$ and $C=\frac{P(\alpha)+Q(\beta)}{D}$ %$A=P(x), B=Q(y)$ and $C=P(x)+Q(y)$ 
yields
\begin{eqnarray*}
\mathrm{max}\{|\frac{P(\alpha)}{D}|,|\frac{Q(\beta)}{D}|, |\frac{C}{D}|\}< L(\epsilon)N(n_1!\cdots n_r!)^{1+\epsilon}< L(\epsilon)(4^{1+\epsilon})^{n_1s}\cdots (4^{1+\epsilon})^{n_rs},
\end{eqnarray*}
where $L(\epsilon)=N(A_1\cdots A_r)\cdot K(\epsilon)$. Here $K(\epsilon)$ denotes the constant obtaind directly from the ABC-conjecture. Since there are only finitely many $(n_1,...,n_r)$ satisfying
\begin{center}
 $n_1!A_1^{n_1}\cdots n_r!A_r^{n_r}=P(\alpha)^2Q(\beta)+Q(\beta)^2P(\alpha)$
\end{center} 
we cannot have infinitely many $\alpha$ and $\beta$ such that $A$, $B$ and $C$ grow constantly. So either we have finitely many $\alpha$ and $\beta$ (in this case we are done) or we have infinitely many $\alpha$ and $\beta$ such that $A,B$ and $C$ do not grow constantly. Notice that there must be infinitely many $\alpha$ and infinitely many $\beta$. Suppose for a moment there were only finitely many, say, $\alpha$ but infinitely many $\beta$. But then $B$ must grow constantly which is excluded. So let us assume that there are infinitely many $\alpha$ and $\beta$ satisfying 
\begin{center}
	$n_1!A_1^{n_1}\cdots n_r!A_r^{n_r}=P(\alpha)^2Q(\beta)+Q(\beta)^2P(\alpha)$.
\end{center} 
Because of 
\begin{eqnarray*}
	\mathrm{max}\{|\frac{P(\alpha)}{D}|,|\frac{Q(\beta)}{D}|, |\frac{C}{D}|\}< L(\epsilon)(4^{1+\epsilon})^{n_1s}\cdots (4^{1+\epsilon})^{n_rs}.
\end{eqnarray*}
this can happen only if there exist constants $s,r$ such that for infinitely many pairs $\alpha, \beta$ we have $sD=P(\alpha)$ and $rD=Q(\beta)$ and hence $Q(\beta)=\frac{r}{s}P(\alpha)$. Such solutions $(n_1,...n_r,\alpha,\beta)$ therefore satisfy
\begin{center}
	$n_1!A_1^{n_1}\cdots n_r!A_r^{n_r}=(\frac{r}{s}+\frac{r^2}{s^2})P(\alpha)^3$.
\end{center}
Now consider the Polynomial $P'(x):=(\frac{r}{s}+\frac{r^2}{s^2})P(x)^3$. According to Theorem 1.2, the diophantine equation
\begin{center}
	$P'(x)=n_1!A_1^{n_1}\cdots n_r!A_r^{n_r}$ 
\end{center}
has only finitely many integer solutions, thus contradicting our assumtion of having infinitely many $\alpha$. We may reverse the role of $\alpha$ and $\beta$. This completes the proof.
 %we can argue by using Thomae's function describing the greatest common divisor to see that there are also only finitely many $x$ and $y$. Using Thomea's funktion actually needs the fact that $x,y$ have to be coprime. This completes the proof.
\begin{exam2}
\textnormal{The ABC-conjecture implies that $n!=x^2y+y^2x$ and $n!=y^2x-x^2y$ have finitely many integer solutions $(n,x,y)$ with $x$ and $y$ being coprime. This does not follow directly from Theorem 1.8, but one can go through the proof to see that for $P(x)=x$ and $Q(y)=y$ the numbers $P(\alpha)$ and $Q(\beta)$ are allways coprime for coprime $\alpha$ and $\beta$. Therefore, there is no need to divide by the greatest common divisor and we have
\begin{eqnarray*}
	\mathrm{max}\{|P(\alpha)|,|Q(\beta)|, |C|\}< L(\epsilon)(4^{1+\epsilon})^{n_1s}\cdots (4^{1+\epsilon})^{n_rs}.
\end{eqnarray*}
This inequality implies that there can be only finitely many coprime $x$ and $y$ satisfying the equation $n!=x^2y+y^2x$. The same is true for $n!=x^2y-y^2x$ or for $(2n)!=x^2y+y^2x$.}
 
\end{exam2}
\section{Proof of Corollary 1.9}
If $a_d+\cdots +a_0>0$ then we can set $x=y$. The diophantine equation simplyfies to $(a_d+\cdots +a_0)\cdot x^d=bn_1!A_1^{n_1}\cdots n_r!A_r^{n_r}$. Theorem 1.1 implies that we can find infinitely many solutions. Now let $a_d>0$. Set $x=sy$ for some big enough positive integer $s$. Then the equation becomes
\begin{center}
	$(a_ds^d+a_{d-1}s^{d-1}+\cdots +a_0)\cdot y^d=bn_1!A_1^{n_1}\cdots n_r!A_r^{n_r}$.
\end{center}
Since we can choose $s$ big enough such that $(a_ds^d+a_{d-1}s^{d-1}+\cdots +a_0)>0$, we conclude with Theorem 1.1 that there are infinitely many integer solutions. The same argument applies to the case $a_0>0$ except that we take $y=sx$ for a big enough positive $s$.
\section{Proof of Theorem 1.12}
Notice that $f(x,y)=x^sy^s(x^s+y^s)$. For a solution $(z_1,...,z_m, n_1,...,n_l,\alpha,\beta)$ with coprime $\alpha$ and $\beta$, the numbers $\alpha^s$, $\beta^s$ and $\alpha^s+\beta^s$ are mutually coprime. We can use weak form of Szpiro's-conjecture and argue as in the proof of Theorem 1.8 to obtain
\begin{eqnarray*}
|bp_1^{z_1}\cdots p_m^{z_m}\cdot n_1!A_1^{n_1}\cdots n_lA_l^{n_l}!|=|f(\alpha,\beta)|< C\cdot p_1\cdots p_m\cdot (4^s)^{n_1}\cdots (4^s)^{n_l},
\end{eqnarray*}
where $C=N(b\cdot A_1\cdots A_l)$. Since 
\begin{eqnarray*}
	4^r(\frac{n_1}{e})^{n_1}\cdots (\frac{n_l}{e})^{n_l}\leq n_1!\cdots n_l!
\end{eqnarray*}
and since 
\begin{center}
	$n_1!\cdots n_l!< |bp_1^{z_1}\cdots p_m^{z_m}\cdot n_1!A_1^{n_1}\cdots n_l!A_l^{n_l}|$
\end{center}
we have 
\begin{eqnarray*}
	 (\frac{n_1}{e})^{n_1}\cdots (\frac{n_l}{e})^{n_l}< C\cdot p_1\cdots p_m\cdot(4^{sn_1-1})\cdots(4^{sn_l-1}). 
\end{eqnarray*}

%\begin{eqnarray*}
%(\frac{n_1}{e})^{n_1}\cdots (\frac{n_r}{e})^{n_r}< C\cdot p_1\cdots p_m\cdot (4^{s\cdot n_1-1})\cdots (4^{s\cdot n_r-1}). 
%\end{eqnarray*}
It is an exercise to conclude that the later inequality is satisfied only for finitely many $(n_1,...,n_l)$. Now applying Thue-Mahler (see \cite{ST}, Theorem 7.2) yields that the equation
\begin{eqnarray*}
bp_1^{z_1}\cdots p_m^{z_m}\cdot n_1!A_1^{n_1}\cdots n_l!A_l^{n_l}=f(x,y)
\end{eqnarray*}
has only finitely many integer solutions $(z_1,...,z_m,n_1,...,n_l,x,y)$ with $\mathrm{gcd}(x,y)=1$ and $z_i\geq 0$.
%\section{Proof of Proposition 1.5}
\begin{exam2}
	\textnormal{Consider the diophantine equation $x^2y\pm y^2x=7^mn!$. Notice that $x^2y\pm y^2x=xy(x\pm y)$ satisfies the assumptions of Theorem 1.11. Then the weak form of Szpiro's-conjecture implies that the equation has finitely many integer solutions $(m,n,x,y)$ with $m,n>0$ and $\mathrm{gcd}(x,y)=1$. Another example is given by $x^4y^2\pm y^4x^2=3^m m!n!$. One can easily verify that that $x^4y^2- y^4x^2$ satisfies the assumptions of Theorem 1.11.}
\end{exam2}

%\section{Proof of Proposition 1.5}
%As in \cite{T}, Section 8 after multiplication by a suitable integer one gets
%\begin{eqnarray}
%	ay^{n!}=(q_1x-r_1)(q_2x-r_2)(a_{m-2}x^{m-2}+...+a_0),
%\end{eqnarray}
%where $a,q_1,r_1,q_2,r_2,a_{m-2},...,a_0$ are integers with $\frac{r_1}{q_1}$, $\frac{r_2}{q_2}$ different, $\mathrm{gcd}(q_1,r_1)=\mathrm{gcd}(q_2,r_2)=1$ and $\mathrm{gcd}(a,a_{m-2},...,a_0)=1$. For a solution $(x_0,y_0)$ of (16) we obtain integers $a_1,y_1$ with $a_1|a$ and $y_1|y_0$ such that
%\begin{eqnarray*}
%	d_1(q_1x_0-r_1)=a_1e_1y_1^{n!}
%\end{eqnarray*}
%and 
%\begin{eqnarray*}
%	d_2(q_2x_0-r_2)=a_2e_2y_2^{n!}.
%\end{eqnarray*}
%This gives 
%\begin{eqnarray*}
%	C_1y_1^{n!}-C_2y_2^{n!}=C_3
%\end{eqnarray*}
%where $C_1,C_2$ and $C_3$ are non-zero integers bounded from above by a constant depending only on $P(x)$ (see \cite{T}). With a result of Baker \cite{BA}, we finally obtain (see \cite{T})
%\begin{eqnarray*}
%	n!\mathrm{log}(y_2)\leq C_4\mathrm{log}(n!)\mathrm{log}(y_2)
%\end{eqnarray*}
%for a suitable constant $C_4$ and thus
%\begin{eqnarray*}
%	\frac{n!}{\mathrm{log}(n!)}\leq C_5.
%\end{eqnarray*}
%Hence $n$ is bounded by a constant depending only on $P(x)$.

\vspace{0.5cm}
\noindent
{\tiny HOCHSCHULE FRESENIUS UNIVERSITY OF APPLIED SCIENCES 40476 D\"USSELDORF, GERMANY.}
E-mail adress: sasa.novakovic@hs-fresenius.de\\
\noindent
{\tiny MATHEMATISCHES INSTITUT, HEINRICH--HEINE--UNIVERSIT\"AT 40225 D\"USSELDORF, GERMANY.}
E-mail adress: novakovic@math.uni-duesseldorf.de


\begin{thebibliography}{999}
%\bibitem{BA} A. Baker: A sharpening of the bounds for linear forms in logarithms. Acta Arith. 21 (1972), 117-129.
\bibitem{BH} M. Bhargava: P-orderings and polynimial functions on arbitrary subsets of Dedekind rings. J. Reine Angew. Math. 490 (1997), 101-127.
\bibitem{BO} D. Berend and C.F. Osgood: On the equation $P(x)=n!$ and a question of Erd\"os. J. Number Theory. 42 (1992), 189-193.
\bibitem{BH1} D. Berend and J. E. Harmse: On polynomial-factorial diophantine equations. Trans. Amer. Math. Soc. 358 (2006), 1741-1779.
\bibitem{BG} B.C. Berndt and W.F. Galeway: On the Brocard-Ramanujan diophantine equation $n!+1=m^2$. The Ramanujan J. 4 (2000), 41-42.
\bibitem{BR} H. Brocard: Question 1532. Nouv. Corresp. Math. 2 (1876); Nouv. Ann. Math. 4 (1885), 391.
\bibitem{EO} P. Erd\"os and R. Obl\'ath: \"Uber diophantische Gleichungen der Form $n!=x^p \pm y^p$ und $n!\pm m!= x^p$. Acta Szeged. 8 (1937), 241-255.
\bibitem{BMS} Y. Bugeaud, M. Mignotte and S. Siksek: Classical and modular aproaches to exponential and diophantine equations II. The Lebesgue--Nagel equation. Compos. Math. 142 (2006), 31-62.
\bibitem{BPZ} HM. Bui, K. Pratt and A. Zaharescu: Power savings for counting solutions to polynomial-factorial equations. 	arXiv:2204.08423. 
\bibitem{DAB} A. Dabrowski: On the equation $n!+A=y^2$. Nieuw Arch. Wisk. 14 (1996), 321-324.
\bibitem{DA} A. Dabrowski: On the Brocard--Ramanujan problem and generalizations. Coll. Mathe. 126 (2012), 105-110.
\bibitem{DMU} A. Dabrowski and M. Ulas: Variations on the Brocard--Ramanujan equation. J. Number Theory 133 (2013), 1168-1185.
\bibitem{EG} A. Epstein and J. Glickman (2020), https://github.com/jhg023/brocard
\bibitem{KL} O. Kihel and F. Luca: Variants of the Brocard--Ramanujan equation. J. Th\'eor. Nombres Bordeaux 20 (2008), 353-363.
\bibitem{LA} S. Lang: Old and new conjectured diophantine inequalities. Bull. Amer. Math. Soc. 23 (1990), 37-75.
\bibitem{L} F. Luca: The diophantine equation $P(x)=n!$ and a result of M. Overholt. Glasnik Matemati\'cki 37 (2002), 269-273.
\bibitem{MA} R. Matson: Brocard's problem 4th solution search utilizing quadratic residues. Unsolved Problems in Number Theory, Logic and Cryptography (2017), available at http://unsolvedproblems.org/S99.pdf
\bibitem{O} M. Overholt: The diophantine equation $n!+1=m^2$. Bull. London. Math. Soc. 42 (1993), 104.
\bibitem{PS} R.M. Pollack and H.N. Shapiro: The next to last case of a factorial diophantine equation. Comm. Pure Appl. Math. 25 (1973), 313-325. 
\bibitem{RA} S. Ramanujan: Question 469. J. Indian Math. Soc. 5 (1913), 59.
\bibitem{SH} T.N. Shorey: Diophantine approximations, diophantine equations, transcendence and applications. Indian Jour. of Pure and Applied Math. 37 (2006), 9-39.
\bibitem{ST} T.N. Shorey and R. Tijdeman: Exponential diophantine equations, Cambridge Tracts in Math. 87, Cambridge University Press (1986).
\bibitem{SS} S. Siksek: Diophantine equations after Fermat's last theorem.  J. Th\'eor. Nombres Bordeaux 21 (2009), 425-436.
\bibitem{WT} W. Takeda: On the finiteness of solutions for polynomial-factorial diophantine equations. Forum Math. 33 (2021), 361-374.
%\bibitem{S} C.L. Siegel: Aproximation algebraischer Zahlen. Math. Zeit. 10 (1921), 173-213.
%\bibitem{T} R. Tijdeman: Applications of the Gel'fond--Baker method to rational number theory. Colloquia Math. Soc. J\'anos Bolyai (13), Topics in number theory (1974), 339-416.
\bibitem{MU} M. Ulas: Some observations on the diophantine equation $y^2=x!+A$ and related results. Bull. Aust. Math. Soc. 86 (2012), 377-388.

\bibitem{MUT} M. Ulas: Some experiments with Ramanujan--Nagell type diophantine equations. Glasnik Matemati\'cki 49 (2014), 287-302.
\bibitem{TY} T. Yamada: A generalization of the Ramanujan--Nagell equation. Glasgow Math. J. 61 (2019), 535-544.


\end{thebibliography}
\end{document}